# Cross-benchmarking for performance evaluation: looking across best practices of different peer groups using DEA

*Nuria Ramón, José L. Ruiz and Inmaculada Sirvent*

*Centro de Investigación Operativa. Universidad Miguel Hernández. Avd. de la Universidad, s/n 03202-Elche (Alicante), SPAIN*

*September, 2019*

**Abstract**

In benchmarking, organizations look outward to examine others' performance in their industry or sector. Often, they can learn from the best practices of some of them and improve. In order to develop this idea within the framework of Data Envelopment Analysis (DEA), this paper extends the common benchmarking framework proposed in Ruiz and Sirvent (2016) to an approach based on the benchmarking of decision making units (DMUs) against several reference sets. We refer to this approach as *cross-benchmarking*. First, we design a procedure aimed at making a selection of reference sets (as defined in DEA), which establish the common framework for the benchmarking. Next, benchmarking models are formulated which allow us to set the closest targets relative to the reference sets selected. The availability of a wider spectrum of targets may offer managers the possibility of choosing among alternative ways for improvements, taking into account what can be learned from the best practices of different peer groups. Thus, cross-benchmarking is a flexible tool that can support a process of future planning while considering different managerial implications.

*Keywords*: Performance Evaluation; Data Envelopment Analysis; Benchmarking; Target Setting.

## 1. Introduction

Data Envelopment Analysis (DEA) (Charnes et al. (1978)) is widely used as a benchmarking tool for improving performance of organizations in an industry or sector. In DEA, the performance of decision making units (DMUs) is evaluated against targets set out on the efficient frontier of the production possibility set (PPS), which can be seen as a best practice frontier in the circumstance of benchmarking (see Cook et al. (2014) for discussions). Comparisons between actual performances and targets may show the DMUs the way for improvement. Some recent papers dealing with DEA and benchmarking include Adler et al. (2013), which uses network DEA, Zanella et al. (2013), which uses a model based on a novel specification of weight restrictions, Dai and Kuosmanen (2014), which combines DEA with clustering methods, Yang et al. (2015), which uses DEA to create a dynamic benchmarking system, Gouveia et al. (2015), which combines DEA and multi criteria decision analysis (MCDA), Daraio and Simar (2016), which deals with benchmarking and directional distances, and Ghahraman and Prior (2016) and Lozano and Calzada-Infante (2018), which propose stepwise benchmarking approaches.

As stated in Thanassoulis et al. (2008), non-radial DEA models are the appropriate instruments for target setting and benchmarking, because they ensure that the identified targets lie on the Pareto-efficient subset of the frontier. Within that family of models, there is already a large and growing body of research on setting the closest targets. Closest targets minimize the gap between actual performances and best practices, thus defining a plan for improvements that requires as little effort as possible from the DMUs for their achievement. Aparicio et al. (2007) developed a mixed integer linear programming model that minimizes the distance to the strong efficient frontier of the PPS. This model allows therefore to setting the closest targets. Since then, the ideas in that paper have been extended to deal with different issues related to the benchmarking: in Ramón et al. (2016), which develops DEA benchmarking models with weight restrictions, in Ruiz and Sirvent (2016), which proposes a common framework for the benchmarking, in Aparicio et al. (2017), which deals with oriented models, in Cook et al. (2017), which develops an approach for the benchmarking of DMUs classified in groups of units that experience similar circumstances, in Ramón et al. (2018), which proposes a stepwise benchmarking, in Ruiz and Sirvent (2019), which incorporates information on goals into the benchmarking, and in Cook et al. (2019), which deals with the benchmarking within the context of pay-for-performance incentive plans. See also Portela et al. (2003) and Tone (2010) for other DEA papers dealing with closest targets.

For performance evaluation, we take into consideration issues like whether or not allowing for individual circumstances of the DMUs (see Roll et al. (1991) for discussions). In the affirmative case, DEA can be used, because it provides an evaluation based on DMU-specific input/output weights. In the negative one, we can use a common set of weights (CSW). This is the traditional approach for efficiency measurement in Economics and Engineering. In this latter case, we may alternatively use cross-efficiency evaluation (Sexton et al. (1986) and Doyle and Green (1994)), which also provides a common evaluation framework, in the sense that DMUs are assessed with several profiles of weights, but these profiles are the same in the evaluation of all of them. A similar reasoning can be followed if target setting and benchmarking is the purpose of the analysis instead of efficiency measurement. We can use DEA, which means to evaluate performance against targets determined by one DMU-specific reference set, or we can use the approach in Ruiz and Sirvent (2016), which implies using one reference set that is common to all units. As an alternative, in the present paper we propose to extend that 2016 study to a benchmarking approach that seeks to evaluate performance against several reference sets. Thus, more flexibility is allowed in the selection of reference sets, but within a common benchmarking framework because target setting is carried out by using the same reference sets for all the DMUs. We refer to this approach as *cross-benchmarking*[1]. The availability of several reference sets is particularly desirable when benchmarking is the purpose of analysis. Considering various reference sets makes it possible to evaluate each unit against a wider spectrum of targets, which may offer different insights on managerial implications. Thus, managers have the possibility of choosing among different alternatives the way for improvement, taking into account what can be learned from the best practices of a broader range of peer groups. Obviously, cross-benchmarking also somehow extends the idea behind cross-efficiency evaluation to the context of benchmarking. In cross-efficiency evaluation, every DMU is evaluated from the perspective of all of the others by using their input/output DEA weights. This makes it possible a peer-evaluation of DMUs, as opposed to the DEA self-evaluation, in which each unit is assessed only with its own weights (those which show them in their best possible light). In the same way, cross-benchmarking provides an evaluation of performance by looking across the different best practices of several reference sets.

One of the key issues that needs to be addressed in cross-benchmarking is the selection of the reference sets the performance of the DMUs is evaluated against. We note that the efficient DMUs in a given reference set span an efficient face of the strong efficient frontier of the PPS, where the DMUs are projected on for performance evaluation. Some studies have shown that the DEA strong efficient frontier in practice consists very often of a large number of efficient facets, so there are in principle many reference sets that can be considered for the benchmarking[2,3]. This is why we develop here a procedure aimed at making a representative selection of reference sets, in a certain sense. Specifically, the proposed procedure seeks to make a selection that ensures no further improvements through considering more reference sets, in terms of reducing the gap between actual performances and best practices. Once the selection of reference sets is made, the benchmarking models we develop follow a criterion of minimization of the distance to the frontier for target setting, so they set the closest targets relative to each of them. Setting the closest targets is a suitable approach for the benchmarking in absence of information on preferences because, as said before, they are the ones that require less effort for improvement. This is therefore a data-driven approach for target setting based on the closeness between actual inputs/outputs and targets (in line with others mentioned above), which has now been developed allowing for several reference sets within a common benchmarking framework.

---

[1] As far as we know, the term "cross-benchmarking" within the DEA literature has only been used in Mousavi and Ouenniche (2018), which follows an approach based on context-dependent DEA (Seiford and Zhu, 2003), in contrast to the one proposed here which uses the efficient frontier determined by the whole sample of DMUs.
[2] We note that efficient faces are not necessarily maximal faces. Maximal efficient faces are usually called facets.
[3] There exists a body of literature dealing with the identification of the facets of the DEA efficient frontier (see Olesen and Petersen (1996, 2003, 2015), Fukuyama and Sekitani (2012), Jahanshahloo et al. (2007) and Davtalab-Olyaie et al. (2014)).

Finally, it should be noted that, although the proposed approach is developed within the context of the DEA methodology, we here move away from the notion of technical efficiency typically used in standard DEA. In other words, in the cross-benchmarking dominance does not prevail, as a result of setting targets from projections on to a given selection (common to all units) of efficient faces of the DEA strong efficient frontier. This means that this approach considers the possibility that the DMUs can learn from benchmarks other than those that are better in all of the dimensions of performance, provided they suggest ways of improvement that can be realistically implemented (according to managers' opinion). As a result, the targets set may suggest improvements through reallocations between inputs and/or outputs, which managers might consider in future planning. Something similar occurs with the DEA models that include weight restrictions, wherein we actually move from technical efficiency to a kind of overall efficiency (see Thanassoulis et al. (2004)).

The paper is organized as follows: in section 2, we design a selection algorithm of reference sets to be used in the cross-benchmarking. In Section 3 we develop models for the benchmarking/target setting against the reference sets previously selected. Section 4 includes an empirical illustration. Last section concludes.

## 2. Selecting reference sets for cross-benchmarking

Throughout the paper, it is supposed that we evaluate the performance of $n$ DMUs that use $m$ inputs to produce $s$ outputs. These are denoted by $(X_j, Y_j)$, $j=1,...,n$, where $X_j = (x_{1j},...,x_{mj})' > 0_m$, $j=1,...,n$, and $Y_j = (y_{1j},...,y_{sj})' > 0_s$, $j=1,...,n$. For the benchmarking, we assume a constant returns to scale (CRS) technology (Charnes et al., 1978). Thus, the PPS, $T = \{(X,Y)/X \text{ can produce } Y\}$, can be described as $T = \left\{ (X,Y) \middle/ X \geq \sum_{j=1}^{n} \lambda_j X_j, Y \leq \sum_{j=1}^{n} \lambda_j Y_j, \lambda_j \geq 0 \right\}$.

In this section, we design a sequential procedure aimed at making a selection of reference sets to be used in the cross-benchmarking of all the DMUs. This includes, in particular, the efficient DMUs, which means that cross-benchmarking contemplates the possibility that efficient DMUs may learn from other efficient units (in that sense, there is a parallelism with cross-efficiency evaluation, wherein efficient DMUs are evaluated with the DEA weights of other units). Target setting will be carried out as the result of projecting the DMUs on to the different faces of the efficient frontier spanned by the members of the reference sets that have been selected. The idea is to add a new reference set at each step to the ones selected in prior steps, by moving from the common benchmarking (all the DMUs are evaluated against the same reference set) to something akin to a self-benchmarking, in which every DMU can eventually choose the most favourable reference set to evaluate against. Specifically, in the first step, $R_1$ is selected as the common reference set that allows us to set the closest targets on to the face spanned by its members, by minimizing globally differences between actual performances (data) and targets (the projections onto such face). That is, $R_1$ provides the best evaluation of the DMUs (globally) within a common benchmarking framework. The addition of a new reference set offers the possibility of finding closer targets for some DMUs. Thus, in step 2, DMUs are allowed to evaluate against either $R_1$ or a new reference set $R_2$, which is selected as the one that leads to the closest targets globally on to the two faces spanned by the corresponding DMUs in these two sets. And so on. The procedure continues to consider a new reference set at each step, which is identified by minimizing the distances between actual input/outputs and targets on to the faces spanned by the reference sets selected, until no DMU can find closer targets through considering more reference sets.

## 2.1. Step 1: Selecting a common reference set

The following model provides the reference set $R_1$, which can be seen as a common reference set for the benchmarking of all the DMUs

$$\text{Min} \quad D_1 = \sum_{j=1}^{n} d_j^1$$

s.t.:

$$\sum_{k \in E} \lambda_{kj}^1 X_k = \hat{X}_j^1 \qquad j=1,...,n \qquad (1.1)$$

$$\sum_{k \in E} \lambda_{kj}^1 Y_k = \hat{Y}_j^1 \qquad j=1,...,n \qquad (1.2)$$

$$d_j^1 = \left\| (X_j, Y_j) - (\hat{X}_j^1, \hat{Y}_j^1) \right\|_1^\omega \qquad j=1,...,n \qquad (1.3)$$

$$-V'X_k + U'Y_k + b_k = 0 \qquad k \in E \qquad (1.4)$$

$$V \geq 1_m \qquad (1.5)$$

$$U \geq 1_s \qquad (1.6)$$

$$\lambda_k^1 = \sum_{j=1}^{n} \lambda_{kj}^1 \qquad k \in E \qquad (1.7)$$

$$\lambda_k^1 b_k = 0 \qquad k \in E \qquad (1.8)$$

$$b_k \geq 0, \lambda_k^1 \geq 0 \qquad k \in E$$

$$\lambda_{kj}^1 \geq 0 \qquad k \in E, j=1,...,n \qquad (1)$$

$$d_j^1 \geq 0, \hat{X}_j^1 \geq 0_m, \hat{Y}_j^1 \geq 0_s \qquad j=1,...,n$$

where $\left\| (X_j, Y_j) - (\hat{X}_j^1, \hat{Y}_j^1) \right\|_1^\omega = \sum_{i=1}^{m} \frac{|x_{ij} - \hat{x}_{ij}^1|}{x_{ij}} + \sum_{r=1}^{s} \frac{|y_{rj} - \hat{y}_{ij}^1|}{y_{rj}}$, $j=1,...,n$, is the weighted $L_1$-distance between actual inputs/outputs and targets, $1_m$ ($0_m$) is the m-vector of 1's (0's) ($1_s$ and $0_s$ can be defined analogously) and E is the set of extreme efficient DMUs of the PPS (see Charnes et al., 1991). As for the variables of the model, they are in most cases the ones typically used in DEA formulations: $\hat{X}_j^1$ ($\hat{Y}_j^1$) is the input (output) vector of the projection on the frontier of DMU$_j$, $\lambda_{kj}^1$ are the intensities associated with such projection, $\lambda_k^1$ is the sum of intensities across all units for a given DMU$_k$ in E, $d_j^1$ represents the deviations between actual inputs/outputs and targets for DMU$_j$, V (U) is the input (output) weight vector and $b_k$ is a variable used to express the classical constraints $-V'X_k + U'Y_k \leq 0$ as equalities.

Some of the constraints of (1) are common to all of the models that set the closest targets, which result from the characterization of the strong efficient frontier of T provided in Aparicio et al. (2007). Through these constraints we consider projections $(\hat{X}_j^1, \hat{Y}_j^1)$, $j=1,...,n$, for every DMU$_j$, j=1,…,n, which lie all on a same face of the strong efficient frontier of T. (1.1)-(1.2) ensure that $(\hat{X}_j^1, \hat{Y}_j^1)$, $j=1,...,n$, belong to T. With (1.4)-(1.6) we allow for all of the supporting hyperplanes of T having non-zero coefficients. Note (1.5)-(1.6) are actually some of the constraints of the dual formulation of the additive model. Through (1.7)-(1.8) the two previous groups of constraints are linked, thus ensuring benchmarks on the Pareto efficient frontier of T. If $\sum_{j=1}^{n} \lambda_{kj}^{1} > 0$ then (1.8) implies $b_k = 0$. Thus, if $DMU_k \in E$ plays an active role as referent in the evaluation of some $DMU_j$, j=1,…,n, then it necessarily belongs to $-V'X + U'Y = 0$. Therefore, the projections $(\hat{X}_j^1, \hat{Y}_j^1)$, j=1,...,n, associated with the feasible solutions of (1), are combinations of $DMU_k$'s $\in E$ which lie all on a same face of the efficient frontier, because these $DMU_k$'s belong all to the same supporting hyperplane of T, which has non-zero coefficients.

In terms of the optimal solutions of (1), $R_1$ is defined as the set of efficient DMUs that participate actively as a referent in the evaluation of some of the units. That is:

**Definition 1.**
$$R_1 = \{DMU_k, k \in E / \lambda_k^{1*} > 0\} = \{DMU_k, k \in E / \lambda_{kj}^{1*} > 0 \text{ for some j, j=1,...,n}\} \quad (2)$$

$R_1$ is also the reference set that allows setting the closest targets (globally), $(\hat{X}_j^{1*}, \hat{Y}_j^{1*})$, to the actual inputs/outputs of all the DMUs. Therefore, model (1) identifies the most similar benchmark performances to the actual performances of the DMUs, that is, those from which to learn with less effort (globally).

It should be noted that this common approach for the benchmarking does not require dominance. Otherwise, the feasibility of (1) could not be ensured, because this model makes all of the DMUs be projected on to the same efficient face of the frontier. For this reason, (1) does not include non-negative slacks but the absolute value of the deviations between actual inputs/outputs and targets, which can be both positive and negative (or zero). Model (1) is very similar to the one proposed in the 2016 Ruiz and Sirvent paper; the main difference is in that we here normalize the deviations by actual inputs/outputs (as it is most commonly done) while in that paper the averages across units of the corresponding inputs/outputs are used. Note that the evaluation of the DMUs against the targets is eventually explained in terms of potential changes to be implemented relative to actual inputs/outputs, so this interpretation is more directly related to the model where the deviations come from.

Remark 1. Non-linear constraints $\lambda_k^1 b_k = 0, k \in E$, can be handled in practice by using Special Ordered Sets (SOS) (Beale and Tomlin, 1970). SOS Type 1 is a set of variables where at most one of them may be nonzero. Therefore, if (1.8) is removed from (1) and each pair of variables $\{\lambda_k^1, b_k\}, k \in E$, is defined as a SOS Type 1, then it is guaranteed that $\lambda_k^1$ and $b_k$ cannot be simultaneously positive for DMU$_k$'s, $k \in E$. CPLEX Optimizer can solve LP problems with SOS (and also LINGO). This treatment of non-linear constraints through SOS variables has already been utilized for solving models like (1) in Ruiz and Sirvent (2016, 2019), Aparicio et al. (2017), Cook et al. (2017) and Cook et al. (2019).

Remark 2. Each absolute value $|x_{ij} - \hat{x}_{ij}^1|$ involved in the weighted L$_1$-distances considered in constraints (1.3) can be replaced by the sum of two non-negative variables $|x_{ij} - \hat{x}_{ij}^1| = p_{ij} + n_{ij}$, where $x_{ij} - \hat{x}_{ij}^1 = p_{ij} - n_{ij}$, in order to make (1) a linear model. The same applies to the absolute values of the differences corresponding to the outputs.

## 2.2. Adding a reference set in step a

As said before, once R$_1$ has been identified, in step 2 a new reference set R$_2$ is selected, minimizing again the differences between actual inputs/outputs and targets, but allowing the DMUs to be evaluated against either R$_1$ or R$_2$. Thus, after running a–1 steps, $a \geq 2$, we have available the reference sets $R_1,...,R_{a-1}$. In addition, for each DMU$_j$, $j = 1,...,n$, we also have both the targets $(\hat{X}_j^{h*}, \hat{Y}_j^{h*})$, $h = 1,...,a-1$, corresponding to these reference sets and the deviations $d_j^{h*}$, $h = 1,...,a-1$, between actual inputs/outputs and each of these targets, measured in terms of the corresponding weighted L$_1$-distance.

In step a, we seek the selection of the reference set R$_a$, which will be identified as the result of minimizing globally the differences between actual inputs/outputs and targets when the DMUs are allowed to be evaluated against any of the reference sets $R_1,...,R_{a-1}$ or $R_a$. In order to do so, we propose to solve the following model

$$\text{Min} \quad D_a = \sum_{j=1}^{n} \delta_j^a$$

s.t.:

$$\sum_{k \in E} \lambda_{kj}^a X_k = \hat{X}_j^a \qquad j=1,...,n \qquad (3.1)$$

$$\sum_{k \in E} \lambda_{kj}^a Y_k = \hat{Y}_j^a \qquad j=1,...,n \qquad (3.2)$$

$$d_j^a = \left\| (X_j, Y_j) - (\hat{X}_j^a, \hat{Y}_j^a) \right\|_1^\omega \qquad j=1,...,n \qquad (3.3)$$

$$d_j^{h*} \leq \delta_j^a + M I_j^h \qquad j=1,...,n,\ h=1,...,a-1 \qquad (3.4)$$

$$d_j^a \leq \delta_j^a + M I_j^a \qquad j=1,...,n \qquad (3.5)$$

$$\sum_{h=1}^{a} I_j^h \leq a-1 \qquad j=1,...,n \qquad (3.6)$$

$$-V'X_k + U'Y_k + b_k = 0 \qquad k \in E \qquad (3.7)$$

$$\lambda_k^a = \sum_{j=1}^{n} \lambda_{kj}^a \qquad k \in E \qquad (3.8)$$

$$V \geq 1_m \qquad (3.9)$$

$$U \geq 1_s \qquad (3.10)$$

$$S_k = \{\lambda_k^a, b_k\} \quad \text{SOS1} \qquad k \in E \qquad (3.11)$$

$$I_j^h \in \{0,1\} \qquad j=1,...,n,\ h=1,...,a \qquad (3)$$

$$d_j^a, \delta_j^a \geq 0 \qquad j=1,...,n$$

$$b_k, \lambda_k^a \geq 0 \qquad k \in E$$

$$\lambda_{kj}^a \geq 0 \qquad k \in E,\ j=1,...,n$$

$$\hat{X}_j^a \geq 0_m,\ \hat{Y}_j^a \geq 0_s \qquad j=1,...,n$$

where M is a big positive quantity.

For every $DMU_j$, $j=1,...,n$, model (3) sets targets $(\hat{X}_j^a, \hat{Y}_j^a)$ associated with the new reference set $R_a$. Variables $d_j^a$, $j=1,...,n$, in (3.3) evaluate the deviations between actual inputs/outputs and each of these new targets, in terms of the corresponding weighted $L_1$-distance. And, in the objective, we minimize, for each $DMU_j$, $j=1,...,n$, the minimum of this deviation, $d_j^a$, and those obtained in the previous steps, $d_j^{h*}$, $h=1,...,a-1$, as shown in the proposition below:

**Proposition 1.** The optimal solution of model (3) satisfy

$$\delta_j^{a*} = \min\{d_j^{h*}, h=1,...,a\},\ j=1,...,n.$$

Proof. For each j, $j=1,...,n,$ constraint $\sum_{h=1}^{a} I_j^h \leq a-1$ in (3) forces at least one of the binary variables $I_j^h$ to be zero. That is, $H_j = \{h / I_j^h = 0, h = 1,...,a\} \neq \emptyset$. Therefore, constraints $d_j^{h*} \leq \delta_j^a + M, h \notin H_j$, $d_j^{h*} \leq \delta_j^a, h \in H_j$, and either $d_j^a \leq \delta_j^a + M$ or $d_j^a \leq \delta_j^a$, depending on whether $a \in H_j$ or not, together with the minimization of $\delta_j^a$ in the objective, lead to $d_j^{h_0*} = \delta_j^{a*} = \min\{d_j^{h*}, h = 1,...,a\}$ for some $h_0 \in H_j$. ∎

The statement in Proposition 1 helps better understand how model (3) works. For each DMU$_j$, $j=1,...,n,$ constraints (3.4)-(3.6) force variable $\delta_j^a$ to be an upper bound of at least one of the differences $d_j^{h*}, h = 1,...,a-1, d_j^a$. The objective minimizes the sum of these upper bounds, so that it actually minimizes for each unit the minimum of the distances between data and each of the a corresponding targets. As a consequence of Proposition 1, $(\hat{X}_j^{h_0*}, \hat{Y}_j^{h_0*})$, as denoted in the proof above, are the closest targets to DMU$_j$ among all those that can be set considering for the benchmarking $R_1,...,R_{a-1}$ and the new reference set selected in step a.

**Proposition 2.** It suffices to set M at

$$M = \max_{\substack{k \in E \\ 1 \leq j \leq n}} \left\| (X_j, Y_j) - (X_k, Y_k) \right\|_1^\omega = \max_{\substack{k \in E \\ 1 \leq j \leq n}} \left( \sum_{i=1}^{m} \frac{|x_{ij} - x_{ik}|}{x_{ij}} + \sum_{r=1}^{s} \frac{|y_{rj} - y_{rk}|}{y_{rj}} \right) \quad (4)$$

in model (3).

Proof. For each j, $j=1,...,n,$ if we define $d_j^2$ at step 2 as

$$d_j^2 = \min \left\{ \left\| (X_j, Y_j) - (\hat{X}_j^2, \hat{Y}_j^2) \right\|_1^\omega / (\hat{X}_j^2, \hat{Y}_j^2) = \sum_{k \in R_2} \lambda_{kj}^2 (X_k, Y_k), \lambda_{kj}^2 \geq 0, k \in R_2 \right\},$$

whichever $R_2$ is, then we have $d_j^2 \leq \max_{k \in R_2} \left\| (X_j, Y_j) - (X_k, Y_k) \right\|_1^\omega \leq \max_{k \in E} \left\| (X_j, Y_j) - (X_k, Y_k) \right\|_1^\omega \leq M$, and so, $d_j^2$ is feasible in (3). Therefore, removing from (3) the feasible solutions $d_j^2$ satisfying $d_j^2 > M$ has no effect on its optimal since, as stated in Proposition 1, the objective minimizes an upper bound of either $d_j^2$ or $d_j^{1*}$.

In addition, $d_j^{1*} \leq \delta_j^2 + MI_j^2$ also holds, given that

$$d_j^{1*} = \min\left\{\left\|(X_j, Y_j) - (\hat{X}_j^1, \hat{Y}_j^1)\right\|_1^\omega \Big/ (\hat{X}_j^1, \hat{Y}_j^1) = \sum_{k \in E}\lambda_{kj}^1(X_k, Y_k), \lambda_{kj}^1 \geq 0, k \in E\right\} \leq \max_{k \in E}\left\|(X_j, Y_j) - (X_k, Y_k)\right\|_1^\omega$$

$\leq M$.

Thus, after running a−1 steps, $a \geq 3$, we have that $d_j^{h*} \leq M, h = 1,...,a-1$. Then, following the same reasoning as in step 2, $d_j^a$, defined as the minimum of the differences between the actual data and the targets determined by all the conical combinations of the efficient DMUs in $R_a$, is feasible and not greater than M. Therefore, removing the feasible solutions $d_j^a > M$ from (3) has no effect on the optimal value of the problem. ∎

It is important to note that, when solving (3), it may happen that $\lambda_{kj}^{a*} > 0$ for a given DMU$_j$ and for some DMU$_k$ in E, while targets $(\hat{X}_j^{a*}, \hat{Y}_j^{a*})$ are not closer to the actual inputs/outputs of DMU$_j$ than the others obtained in previous steps, $(\hat{X}_j^{h*}, \hat{Y}_j^{h*})$, $h = 1,...,a-1$ (it will occur when $a \notin H_j$). This means that considering DMU$_k$ for the benchmarking of DMU$_j$ in step a does not provide any improvement in its evaluation (with respect to the evaluations relative to the targets of previous steps). For this reason, the reference set $R_a$ is defined in terms of the optimal solutions of (3) as follows:

**Definition 2.**
$$R_a = \left\{DMU_k, k \in E / \lambda_{kj}^{a*} > 0 \text{ and } \delta_j^{a-1*} > \delta_j^{a*} = d_j^{a*} \text{ for some j, j=1,...,n}\right\} \quad (5)$$

where $\delta_j^{1*} = d_j^{1*}$, $j = 1,...,n$.

That is, $R_a$ consists exclusively of the efficient DMU$_k$'s that allow us to find closer targets for some unit (with respect to those obtained in previous steps).

Obviously, when the optimal value of model (3) in a given step a equals the one obtained in step a−1, this means that adding a new reference set cannot improve the evaluation of any DMU in terms of finding closer targets. Therefore, the search of reference sets stops, $R_1,...,R_{a-1}$ being the reference sets selected for the cross-benchmarking.

## 2.3. Simplifying model (3)

Model (3) used in step a, $a \geq 2$, can be significantly simplified by using the information collected in the previous step, as shown in the next proposition

**Proposition 3.** Let $R_1, ..., R_{a-1}$ be the reference sets selected in steps $1, ..., a-1$, $a \geq 2$, and set the scalar $\delta_j$ as $\delta_j = \delta_j^{a-1*}$, $j = 1, ..., n$. Then, solving (3) in step a is equivalent to solving

$$\text{Min} \quad D_a = \sum_{j \notin R} \delta_j^a$$

s.t.:

$$\sum_{k \in E} \lambda_{kj}^a X_k = \hat{X}_j^a \quad j \notin R \quad (6.1)$$

$$\sum_{k \in E} \lambda_{kj}^a Y_k = \hat{Y}_j^a \quad j \notin R \quad (6.2)$$

$$d_j^a = \left\| (X_j, Y_j) - (\hat{X}_j^a, \hat{Y}_j^a) \right\|_1^\omega \quad j \notin R \quad (6.3)$$

$$\delta_j \leq \delta_j^a + M I_j \quad j \notin R \quad (6.4)$$

$$d_j^a \leq \delta_j^a + M(1 - I_j) \quad j \notin R \quad (6.5)$$

$$-V' X_k + U' Y_k + b_k = 0 \quad k \in E \quad (6.6)$$

$$\lambda_k^a = \sum_{j=1}^n \lambda_{kj}^a \quad k \in E \quad (6.7)$$

$$V \geq 1_m \quad (6.8)$$

$$U \geq 1_s \quad (6.9)$$

$$S_k = \{\lambda_k^a, b_k\} \quad \text{SOS1} \quad k \in E \quad (6.10)$$

$$I_j \in \{0, 1\} \quad j \notin R$$

$$d_j^a, \delta_j^a \geq 0 \quad j \notin R \quad (6)$$

$$b_k, \lambda_k^a \geq 0 \quad k \in E$$

$$\lambda_{kj}^a \geq 0 \quad k \in E, j \notin R$$

$$\hat{X}_j^a \geq 0_m, \hat{Y}_j^a \geq 0_s \quad j \notin R$$

where $R = \bigcup_{h=1}^{a-1} R_h$.

Proof. Since for each j, $j = 1, ..., n$, $\delta_j$ is defined as $\delta_j = \delta_j^{a-1*} = \min\{d_j^{h*}, h = 1, ..., a-1\}$, then $\delta_j^{a*} = \min\{d_j^{a*}, \delta_j\}$. Therefore, for each j, $j = 1, ..., n$, constraints $d_j^{h*} \leq \delta_j^a + M I_j^h$, $h = 1, ..., a-1$, $d_j^a \leq \delta_j^a + MI_j^a$, and $\sum_{h=1}^{a} I_j^h \leq a-1$ in (3) can be replaced with $\delta_j \leq \delta_j^a + MI_j$ and $d_j^a \leq \delta_j^a + M(1-I_j)$, $I_j$ being a binary variable that takes the value 1 when DMU$_j$'s targets in step a are not closer to this unit than those found in earlier steps.

Moreover, for every DMU$_j$ in R we have $\delta_j = \delta_j^{a-1*} = 0$. Therefore, $\delta_j^{a*} = \min\{d_j^{a*}, 0\} = 0$, and so, for every $j \in R$, constraints $\sum_{k \in E} \lambda_{kj}^a X_k = \hat{X}_j^a$, $\sum_{k \in E} \lambda_{kj}^a Y_k = \hat{Y}_j^a$, $d_j^a = \left\| (X_j, Y_j) - (\hat{X}_j^a, \hat{Y}_j^a) \right\|_1^\omega$, $\delta_j \leq \delta_j^a + MI_j$, $d_j^a \leq \delta_j^a + M(1-I_j)$, can be removed from (3), as well as all the variables associated with $j \in R$. ∎

Model (6) is computationally much simpler than model (3), mainly as a result of reducing notably the number of binary variables that are used. Proof above reveals that 1) it is not necessary to consider the DMUs that belong to reference sets selected in previous steps, and 2) each variable $\delta_j^a$, $j=1,...,n$, can be represented as an upper bound of the scalar $\delta_j$ and the variable $d_j^a$. As a result, model (6) 1) does not include either the variables nor constraints associated with DMUs in R; and, what is more important, 2) only includes one binary variable per DMU to be evaluated (instead of a). In addition, constraints (3.4) disappear as well as (3.6), while (6) only includes in return a new constraint per each unit to be evaluated, in which we have one binary variable.

## 2.4. An algorithm for the selection of reference sets

The procedure we propose for the selection of reference sets to be used in the cross-benchmarking can be implemented applying the following algorithm

1. Identify E and set $R = \varnothing$ and a=1.
2. Solve model (1).
3. Set R$_1$ as in (2), $\delta_j = \delta_j^{1*} = d_j^{1*}$, $j=1,...,n$, D= $D_1^*$, $R = R \cup R_1$ and a=a+1.
4. Solve (6). If $D_a^* = D$, then STOP.
5. Set $R_a$ as in (5), $\delta_j = \delta_j^{a*}$, $j=1,...,n$, D= $D_a^*$, $R = R \cup R_a$, a=a+1, and go to 4.

After running the algorithm, if it took A+1 steps, we will have $R_1,...,R_A$ reference sets to be used in the cross-benchmarking. Obviously, this algorithm finishes in a finite number of steps, as the strong efficient frontier of the PPS consists of a finite number of faces.

## 3. Cross-benchmarking for performance evaluation

Cross-benchmarking, as proposed here, seeks the evaluation of performance of the DMUs against the reference sets selected following the procedure described in the previous section. The evaluation is carried out by comparing actual performances (inputs/outputs) with the targets determined by using each of the reference sets. Targets are determined as projections points on to each of the efficient faces spanned by the DMUs in each of the reference sets. Specifically, we follow an approach based on closest targets, so targets for a given DMU$_j$ associated with the reference set $R_h$, $h=1,...,A$, are determined as the coordinates of the closest projection point on to the efficient face spanned by the DMUs in $R_h$. The following model provides the targets associated with a given reference set $R_h$ for all of the DMUs simultaneously

$$\text{Min} \quad \sum_{j=1}^{n} \left\| (X_j, Y_j) - (\tilde{X}_j^h, \tilde{Y}_j^h) \right\|_1^\omega$$

s.t.:
$$\sum_{k \in R_h} \lambda_{kj}^h X_k = \tilde{X}_j^h \qquad j=1,...,n$$
$$\sum_{k \in R_h} \lambda_{kj}^h Y_k = \tilde{Y}_j^h \qquad j=1,...,n \qquad (7)$$
$$\tilde{X}_j^h \geq 0_m, \tilde{Y}_j^h \geq 0_s \qquad j=1,...,n$$
$$\lambda_{kj}^h \geq 0 \qquad j=1,...,n, k \in R_h$$
$$d_j^h \geq 0 \qquad j=1,...,n$$

In practice, we need to solve model (7) A times, once per reference set. Eventually, we will have A target bundles for each DMU$_j$, which can be found in terms of the optimal solutions of this model as follows

$$\sum_{k \in R_h} \lambda_{kj}^{h*} X_k = \tilde{X}_j^{h*}, \quad h=1,...,A$$
$$\sum_{k \in R_h} \lambda_{kj}^{h*} Y_k = \tilde{Y}_j^{h*}, \quad h=1,...,A \qquad (8)$$

Comparing actual performances (inputs/outputs) and targets that represent best practices of a broader range of peers may offer a wider view on DMUs' performance. In fact, as we shall see later in the empirical illustration, cross-benchmarking often provides different insights about unit's performance across reference sets. Cross-benchmarking may identify units whose actual performance is close to the best practices in almost all of the scenarios of evaluation or, on the contrary, units for which there is a big gap between actual inputs/outputs and targets irrespective of the reference sets they are evaluated against. In some cases, the evaluations suggest the need for improvement in different aspects of performance, while on other occasions, ways for improvements of some inputs and/or outputs can be suggested but in return for the worsening of some others. This all may constitute a valuable source of information in management for decision-making.

## 4. Illustrative example

For purposes of illustration only, in this section we apply the proposed approach to the data set in Coelli et al. (2002), which deals with the evaluation of performance of 28 international airlines during the year 1990. This data set has been widely utilized in the DEA literature, in particular for illustrating the use of methods aimed at both target setting and benchmarking. The airlines are evaluated in terms of two outputs and four inputs. The outputs are: $y_1$=passenger-kilometers flown (PASS) and $y_2$=freight tonne-kilometers flown (CARGO). And the inputs: $x_1$=number of employees (LAB), $x_2$= fuel in millions of gallons (FUEL), $x_3$=other kind of inputs (millions of U.S. dollar equivalent) excluding labour and fuel expenses (MATL) and $x_4$=Capital (CAP), as the sum of maximum takeoff weights of all aircraft flown multiplied by the number of days flown. An initial DEA analysis identifies 9 airlines as efficient: JAL, QANTAS, SAUDIA, SINGAPORE, FINNAIR, LUFTHANSA, SWISSAIR, PORTUGAL and AM.WESTERN.

Regarding the cross-benchmarking, Table 1 reports the results of the algorithm used for the selection of reference sets. We can see that a selection of 7 reference sets is made for the benchmarking after running the corresponding steps of the algorithm. Table 1 specifically shows the airlines included in each of the reference sets. In the last row, the optimal value of the model solved in each step is recorded. In most cases, the reference sets consist of 4 airlines; this happens with R1, R2, R4, R6 and R7. In addition, it is worth highlighting the fact that SINGAPORE belongs to these 5 reference sets. FINNAIR and JAL belong to 4 reference sets, SWISSAIR and AM.WESTERN to 3, while, in contrast, QANTAS and PORTUGAL are only members of R1.

Table 1. Selection of reference sets

For some airlines (which are selected as representative cases), Table 2 records their actual inputs/outputs and the targets obtained by solving model (7) for each of the 7 reference sets. This table also reports the corresponding deviations in percent (between brackets): for the inputs, these are calculated as $\left(x_{ij} - \tilde{x}_{ij}^{h*}\right)/x_{ij}$, i=1,2,3,4, h=1,…,7, and for the outputs as $\left(\tilde{y}_{rj}^{h*} - y_{rj}\right)/y_{rj}$, r=1,2, h=1,…,7. Note that a positive value of the deviations so calculated indicates that the actual input/output is worse than the corresponding target, while a negative value means that the target has been outperformed. Figures 1 to 5 depict graphically the cross-benchmarking of some airlines in terms of such deviations.

Table 2. Actual inputs/outputs and targets by reference set

Figure 1 shows the case of NIPPON. The evaluation against the targets provided by R1, R4 and R6 suggests that others have used less of inputs FUEL, MATL and CAP (by more than 30% in most cases), while maintaining the level of outputs, albeit the number of employees (LAB) is in those cases higher. Moreover, the benchmarking against R2 and R5 reveals that it is possible to operate at a level still lower in some of these inputs and maintaining LAB, although flying less passenger-kilometers (PASS) too. NIPPON managers might find these results useful from the perspective of managerial implications. In particular, looking at others' best practices might be showing how resources can be reallocated (maintaining output generation) or the possible effects of saving resources on passenger-kilometers and/or tonne-kilometers flown. Thus, cross-benchmarking may support a decision making process of planning improvements.

Figure 1. NIPPON: Deviations actual data-targets by reference set

Figure 2 depicts the case of TWA. This figure shows a relatively good performance of this airline in many scenarios. Nevertheless, the evaluation with respect to R3 reveals that it is possible to use less than half of inputs FUEL and MATL, although flying half of passenger-kilometers (PASS) in return. The benchmarking against R7 is somewhat striking, because of the large difference between actual data and target in CARGO. This can be explained by the fact that JAL and LUFTHANSA, which are the two airlines with the highest levels in CARGO, belong to R7. In such a situation, TWA managers may find it unrealistic to emulate the pattern of performance determined by the R7 targets to improve. However, for other airlines flying more tonne-kilometers, like BRITISH, the evaluation against R7 might offer useful information: the green line in Figure 3 now shows that the gap between BRITISH's actual CARGO and that of others having operated in all inputs and outputs at the same level (and even at a lower level in MATL) is not so large. In general, the fact that others have performed in a certain manner does not necessarily mean that a particular organization can realistically emulate all of them. Managers will make the decision on the most appropriate way for improvement by choosing among alternatives.

Figure 2. TWA: Deviations actual data-targets by reference set

Figure 3. BRITISH: Deviations actual data-targets by reference set

Figure 4 depicts the case of AUSTRIA, where we can see that all the points are above the horizontal level 0 line in the benchmarking of this airline against almost all of the reference sets (except R3). Specifically, there is clear room for improvement in inputs LAB and MATL and output PASS. MALAYSIA and GARUDA also show a poor performance (see Table 2).

Figure 4. AUSTRIA: Deviations actual data-targets by reference set

The fact that some of the DEA efficient airlines belong to several reference sets means that they will be given a good evaluation in many scenarios, because actual performances and targets coincide in those cases. This happens to SINGAPORE, which belongs to R1, R2, R4, R6 and R7. JAL also belongs to 4 reference sets. However, the evaluation of this airline with respect to R6 (see Figure 5) suggests that inputs FUEL and CAP could be significantly reduced (by more than 30%), and that such reduction would only entail a loss in CARGO of 10% (the benchmarking against R1 shows a similar pattern). Again, this might be useful information for managers to consider in planning future.

Figure 5. JAL: Deviations actual data-targets by reference set

## 5. Conclusions

Benchmarking involves looking outside a particular organization to examine others' performance in their industry or sector. The ultimate aim is to identify best practices from which to learn and set plans for improvement. The present paper has shown how to implement this core idea through the notion of cross-benchmarking, which is an approach developed within the framework of DEA. Specifically, cross-benchmarking here makes a selection of reference sets, common to all units, which allows us to set a wide spectrum of targets for evaluating DMUs' performance against. Cross-benchmarking is therefore a flexible decision support tool, which allows managers to choose among alternatives when planning improvements, taking into account the best practices of a broad range of peer groups. The results in the empirical application illustrate the fact that managers may set their agenda for improvements by choosing among alternative plans after considering the different implications of reducing input resources and/or increasing output generation, reallocating resources or substituting inputs and outputs.

There are also limitations in cross-benchmarking, which open ways for future research. Some of these limitations are inherited from standard DEA. For instance, there often exists in practice a big gap between actual performances and targets. This raises the need to extend cross-benchmarking for setting more realistically achievable targets, for example by using context-dependent DEA as in Ramón et al. (2018). Moreover, DEA suffers from statistical limitations, as pointed out in Assaf et al. (2011), who use an innovative Bayesian frontier methodology to analyze the performance determinants of retail stores. Cross-benchmarking has been developed here for use when information on preferences is not available. Therefore, it should be extended to allow for incorporating decision maker preferences into the benchmarking (when available). Different approaches can be considered to that end: identifying "model" DMUs (see Brockett et al., 1997), through the specification of a value function (see Halme et al., 1999), adding weight restrictions (see Ramón et al., 2016) or incorporating goal information (see Stewart (2010) and Ruiz and Sirvent (2019)). Likewise, the possibility of dealing with imprecise data should be considered as well. In order to do so, following a fuzzy approach as in León et al. (2003) or using robust optimization as in Toloo and Mensah (2019) are potential alternatives. Finally, we note that cross-benchmarking has been developed assuming a constant returns to scale technology, as has been the case with the existing research on cross-efficiency evaluation (with the exception of Lim and Zhu, 2015). The possible extension of the idea behind the cross-benchmarking to the variable returns to scale case should also be investigated.

**Acknowledgments**

This research has been supported through Grant MTM2016-76530-R (AEI/FEDER, UE).

| Airline | R1 | R2 | R3 | R4 | R5 | R6 | R7 |
| --- | --- | --- | --- | --- | --- | --- | --- |
| JAL |  | ✓ |  | ✓ | ✓ |  | ✓ |
| QANTAS | ✓ |  |  |  |  |  |  |
| SAUDIA |  |  | ✓ |  | ✓ |  |  |
| SINGAPORE | ✓ | ✓ |  | ✓ |  | ✓ | ✓ |
| FINNAIR |  | ✓ |  | ✓ | ✓ | ✓ |  |
| LUFTHANSA |  |  | ✓ |  |  |  | ✓ |
| SWISSAIR |  | ✓ |  |  |  | ✓ | ✓ |
| PORTUGAL | ✓ |  |  |  |  |  |  |
| AM. WEST | ✓ |  |  | ✓ |  | ✓ |  |
| Optimal value | 17.021 | 12.999 | 11.076 | 10.119 | 10.060 | 10.017 | 9.994 |

Table 1. Selection of reference sets

| AIRLINE | INPUTS | | | | OUTPUTS | |
|---|---|---|---|---|---|---|
| | LAB | FUEL | MATL | CAP | PASS | CARGO |
| NIPPON | 12222 | 860 | 2008 | 6074 | 35261 | 614 |
| R1 | 21026.7 (-72%) | 588.9 (32%) | 1250.8 (38%) | 4207.8 (31%) | 35261 (0%) | 614 (0%) |
| R2 | 12222 (0%) | 365 (58%) | 635.1 (68%) | 3744.8 (38%) | 17826.2 (-49%) | 614 (0%) |
| R3 | 19961.5 (-63%) | 189.9 (80%) | 651.2 (68%) | 5509 (9%) | 15325 (-57%) | 614 (0%) |
| R4 | 21460.1 (-76%) | 627.4 (27%) | 1228.6 (39%) | 5082.3 (16%) | 35261 (0%) | 614 (0%) |
| R5 | 12222 (0%) | 365 (58%) | 635.1 (68%) | 3744.8 (38%) | 17826.2 (-49%) | 614 (0%) |
| R6 | 21026.7 (-72%) | 588.9 (32%) | 1250.8 (38%) | 4207.8 (31%) | 35261 (0%) | 614 (0%) |
| R7 | 12222 (0%) | 509.2 (41%) | 1024.5 (49%) | 6074 (0%) | 22637.3 (-36%) | 1508.8(146%) |
| TWA | 35783 | 1118 | 2389 | 8704 | 62345 | 1119 |
| R1 | 40540.4 (-13%) | 1118 (0%) | 2389 (0%) | 7972.6 (8%) | 66986.5 (7%) | 1119 (0%) |
| R2 | 48586 (-36%) | 1081.6 (3%) | 1861.8 (22%) | 8704 (0%) | 58623.1 (-6%) | 1119 (0%) |
| R3 | 35783 (0%) | 343.3 (69%) | 1175.8 (51%) | 9859 (-13%) | 27546.2 (-56%) | 1119 (0%) |
| R4 | 37585.1 (-5%) | 1095.8 (2%) | 2188.8 (8%) | 8704 (0%) | 62345 (0%) | 1119 (0%) |
| R5 | 49466 (-38%) | 1118 (0%) | 1852 (23%) | 9415.4 (-8%) | 58997.8 (-5%) | 1119 (0%) |
| R6 | 40185 (-12%) | 1118 (0%) | 2360.8 (1%) | 7962.5 (9%) | 66901.6 (7%) | 1119 (0%) |
| R7 | 35783 (0%) | 977.1 (13%) | 2389 (0%) | 8704 (0%) | 50840.6 (-19%) | 3317.3(197%) |
| BRITISH | 51802 | 1294 | 4276 | 12161 | 67364 | 2618 |
| R1 | 51802 (0%) | 1294 (0%) | 4149.3 (3%) | 10305.7 (15%) | 81579.5 (21%) | 2618 (0%) |
| R2 | 51802 (0%) | 1325.1 (-2%) | 2539.3 (41%) | 12161 (0%) | 67364 (0%) | 2618 (0%) |
| R3 | 51802 (0%) | 660.2 (49%) | 2164.8 (49%) | 13372.9 (-10%) | 43936.4 (-35%) | 2618 (0%) |
| R4 | 51802 (0%) | 1419.9 (-10%) | 2894.4 (32%) | 12161 (0%) | 78613.3 (17%) | 2618 (0%) |
| R5 | 51802 (0%) | 1298.6 (0%) | 2549.3 (40%) | 16790.3 (-38%) | 67364 (0%) | 2618 (0%) |
| R6 | 51802 (0%) | 1294 (0%) | 2731.1 (36%) | 10463.8 (14%) | 71979.1 (7%) | 2618 (0%) |
| R7 | 51802 (0%) | 1341.1 (-4%) | 3166.0 (26%) | 12161 (0%) | 67364 (0%) | 4482.6 (71%) |
| AUSTRIA | 4067 | 62 | 241 | 587 | 2943 | 65 |
| R1 | 3562 (12%) | 62 (0%) | 241 (0%) | 483.3 (18%) | 4040.6 (37%) | 67.5 (4%) |
| R2 | 2696.7 (34%) | 62 (0%) | 103.1 (57%) | 530.1 (10%) | 3254.8 (11%) | 65 (0%) |
| R3 | 2126.9 (48%) | 20.2 (67%) | 69.4 (71%) | 587 (0%) | 1632.9 (-45%) | 65.4 (1%) |
| R4 | 2696.7 (34%) | 62 (0%) | 103.1 (57%) | 530.1 (10%) | 3254.8 (11%) | 65 (0%) |
| R5 | 3072.3 (25%) | 62 (0%) | 108.9 (55%) | 587 (0%) | 3414.3 (16%) | 65 (0%) |
| R6 | 2912.4 (28%) | 62 (0%) | 105.4 (56%) | 494.2 (16%) | 3315.1 (13%) | 65 (0%) |
| R7 | 2341.1 (42%) | 61.3 (1%) | 137.7 (43%) | 587 (0%) | 2943 (0%) | 199.1 (206%) |
| MALAYSIA | 15156 | 279 | 1246 | 2258 | 12891 | 599 |
| R1 | 14929.5 (2%) | 261.3 (6%) | 1246 (0%) | 2258 (0%) | 17696.5 (37%) | 599 (0%) |
| R2 | 13026.6 (14%) | 279 (0%) | 562.7 (55%) | 2258 (0%) | 14480.3 (12%) | 599 (0%) |
| R3 | 15156 (0%) | 194.6 (30%) | 637.5 (49%) | 3904.5 (-73%) | 12891 (0%) | 774.9 (29%) |
| R4 | 9653.1 (36%) | 279 (0%) | 619.8 (50%) | 2306. 9 (-2%) | 16047.5 (25%) | 599 (0%) |
| R5 | 12880.5 (15%) | 279 (0%) | 458 (63%) | 2258 (0%) | 14918.9 (16%) | 245.3 (-60%) |
| R6 | 13605.3 (10%) | 279 (0%) | 568.9 (54%) | 2161.6 (4%) | 14641.9 (14%) | 599 (0%) |
| R7 | 13174.7 (13%) | 279 (0%) | 668.7 (46%) | 2258 (0%) | 14035.8 (9%) | 957.5 (60%) |
| GARUDA | 10428 | 304 | 3171 | 3305 | 14074 | 539 |
| R1 | 10428 (0%) | 304 (0%) | 745.1 (77%) | 2429 (27%) | 18814.4 (34%) | 539 (0%) |
| R2 | 10428 (0%) | 304 (0%) | 539.7 (83%) | 3079.8 (7%) | 14879.9 (6%) | 539 (0%) |
| R3 | 12481.8 (-20%) | 144.1 (53%) | 479.1 (85%) | 3305 (0%) | 10213.2 (-27%) | 539 (0%) |
| R4 | 10428 (0%) | 315.4 (-4%) | 549.8 (83%) | 3258.7 (1%) | 15355.2 (9%) | 539 (0%) |
| R5 | 10428 (0%) | 304 (0%) | 543 (83%) | 3297.2 (0%) | 14979.1 (6%) | 539 (0%) |
| R6 | 10428 (0%) | 283 (7%) | 597.8 (81%) | 2326.3 (30%) | 16073.1 (14%) | 539 (0%) |
| R7 | 9515.7 (9%) | 304 (0%) | 648.4 (80%) | 3259.7 (1%) | 14074 (0%) | 945.5 (75%) |
| SINGAPORE | 10864 | 523 | 1512 | 4479 | 32404 | 1902 |
| R1 | 10864 (0%) | 523 (0%) | 1512 (0%) | 4479 (0%) | 32404 (0%) | 1902 (0%) |
| R2 | 10864 (0%) | 523 (0%) | 1512 (0%) | 4479 (0%) | 32404 (0%) | 1902 (0%) |
| R3 | 21858.5 (-101%) | 489.7 (6%) | 1512 (0%) | 4479 (0%) | 23789.9 (-27%) | 2395.8 (26%) |
| R4 | 10864 (0%) | 523 (0%) | 1512 (0%) | 4479 (0%) | 32404 (0%) | 1902 (0%) |
| R5 | 10864 (0%) | 523 (0%) | 961.8 (36%) | 6502.9 (-45%) | 23116 (-29%) | 1301 (-32%) |
| R6 | 10864 (0%) | 523 (0%) | 1512 (0%) | 4479 (0%) | 32404 (0%) | 1902 (0%) |
| R7 | 10864 (0%) | 523 (0%) | 1512 (0%) | 4479 (0%) | 32404 (0%) | 1902 (0%) |
| JAL | 21430 | 1351 | 2536 | 17932 | 57290 | 3781 |
| R1 | 21430 (0%) | 1031.7 (24%) | 2982.5 (-18%) | 8835.1 (51%) | 63919.2 (12%) | 3751.8 (-1%) |
| R2 | 21430 (0%) | 1351 (0%) | 2536 (0%) | 17932 (0%) | 57290 (0%) | 3781 (0%) |
| R3 | 32190.1 (-50%) | 762.4 (44%) | 2343.9 (8%) | 6368.2 (65%) | 36062.4 (-37%) | 3781 (0%) |
| R4 | 21430 (0%) | 1351 (0%) | 2536 (0%) | 17932 (0%) | 57290 (0%) | 3781 (0%) |
| R5 | 21430 (0%) | 1351 (0%) | 2536 (0%) | 17932 (0%) | 57290 (0%) | 3781 (0%) |
| R6 | 21430 (0%) | 936 (31%) | 2677.7 (-6%) | 7944.3 (56%) | 57290 (0%) | 3395.5 (-10%) |
| R7 | 21430 (0%) | 1351 (0%) | 2536 (0%) | 17932 (0%) | 57290 (0%) | 3781 (0%) |

Table 2. Actual inputs/outputs and targets by reference set

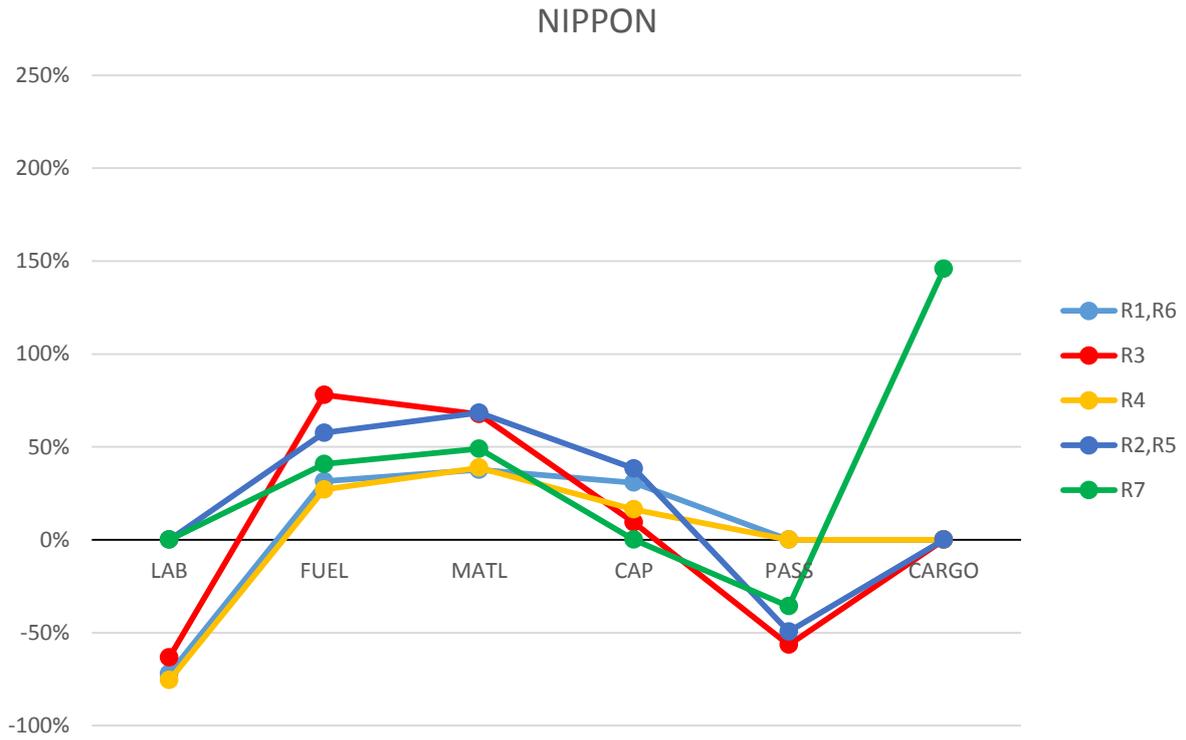

Figure 1. NIPPON: Deviations actual data-targets by reference set

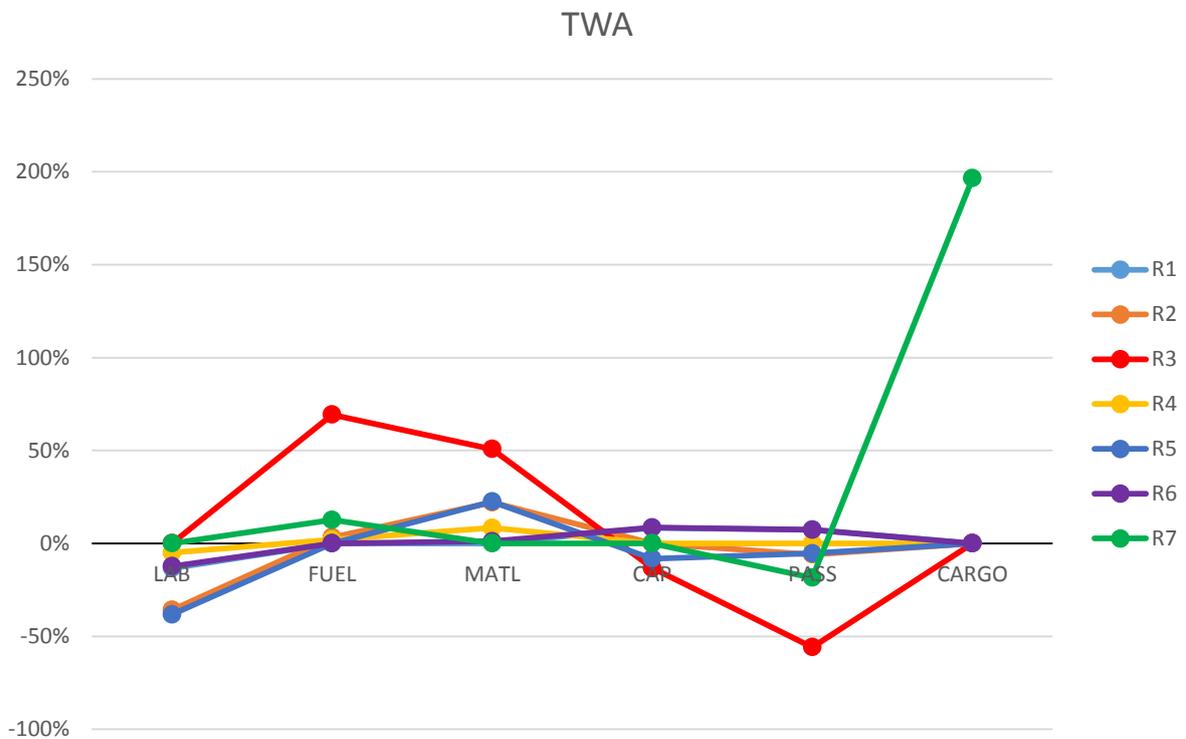

Figure 2. TWA: Deviations actual data-targets by reference set

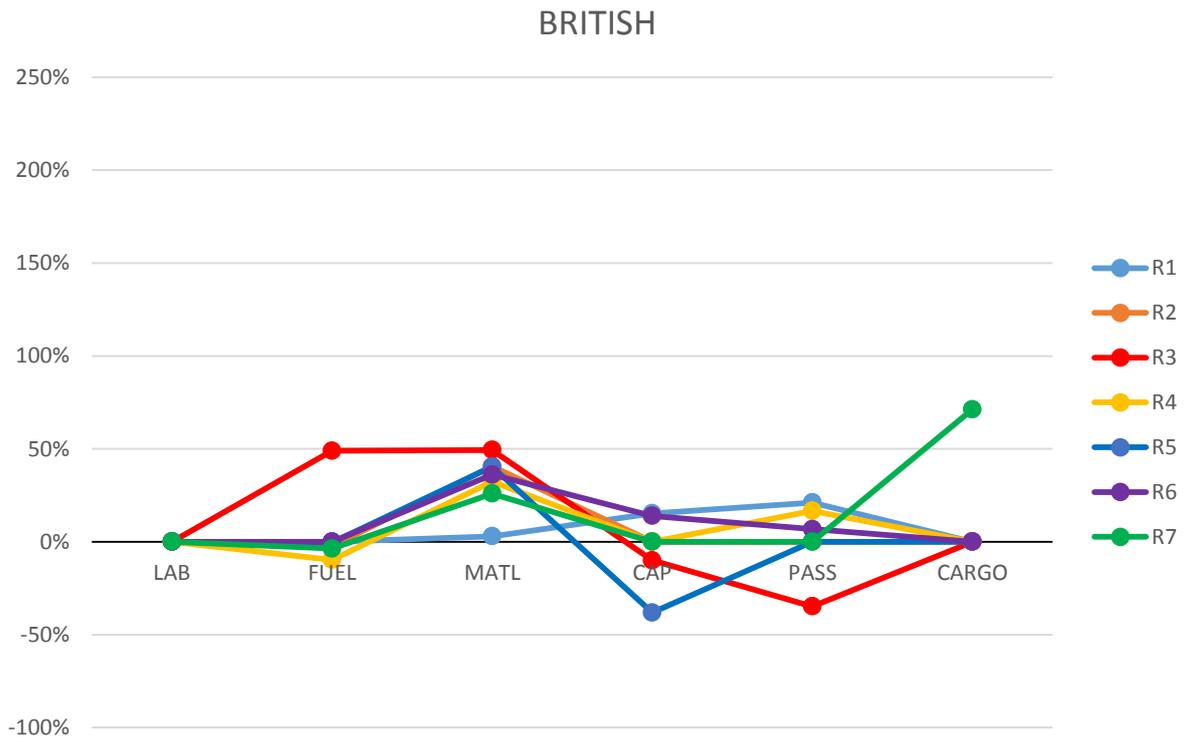

Figure 3. BRITISH: Deviations actual data-targets by reference set

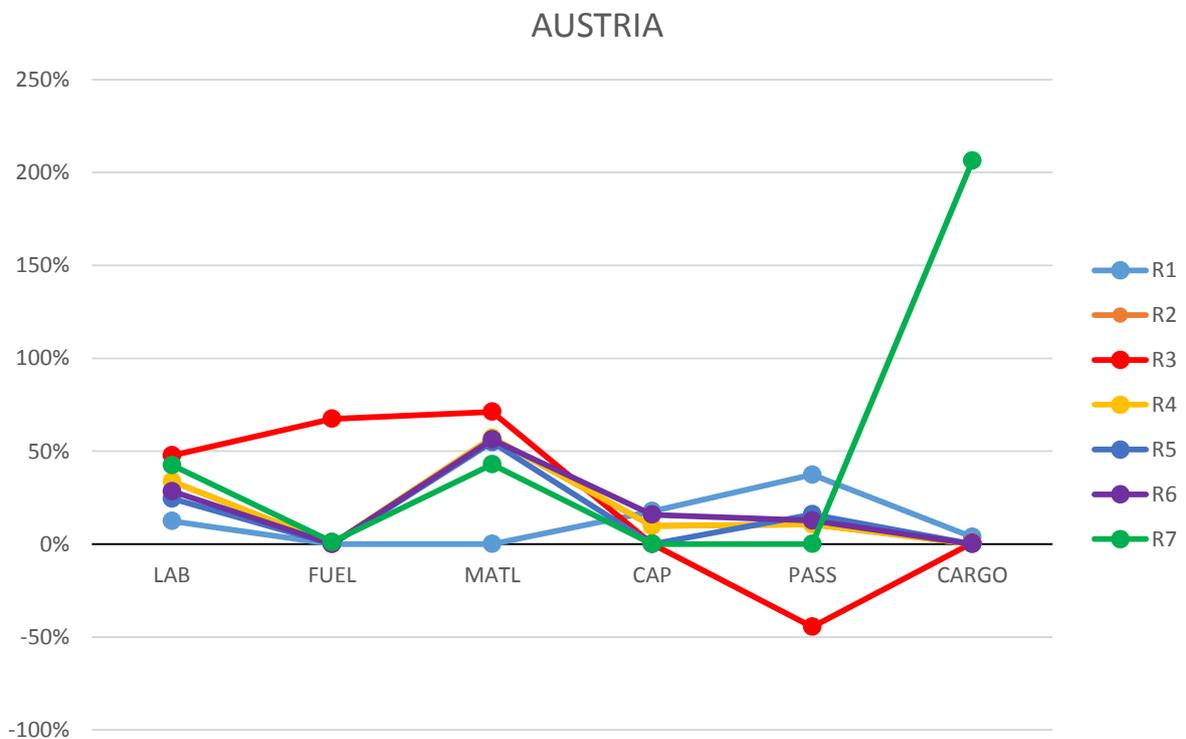

Figure 4. AUSTRIA: Deviations actual data-targets by reference set

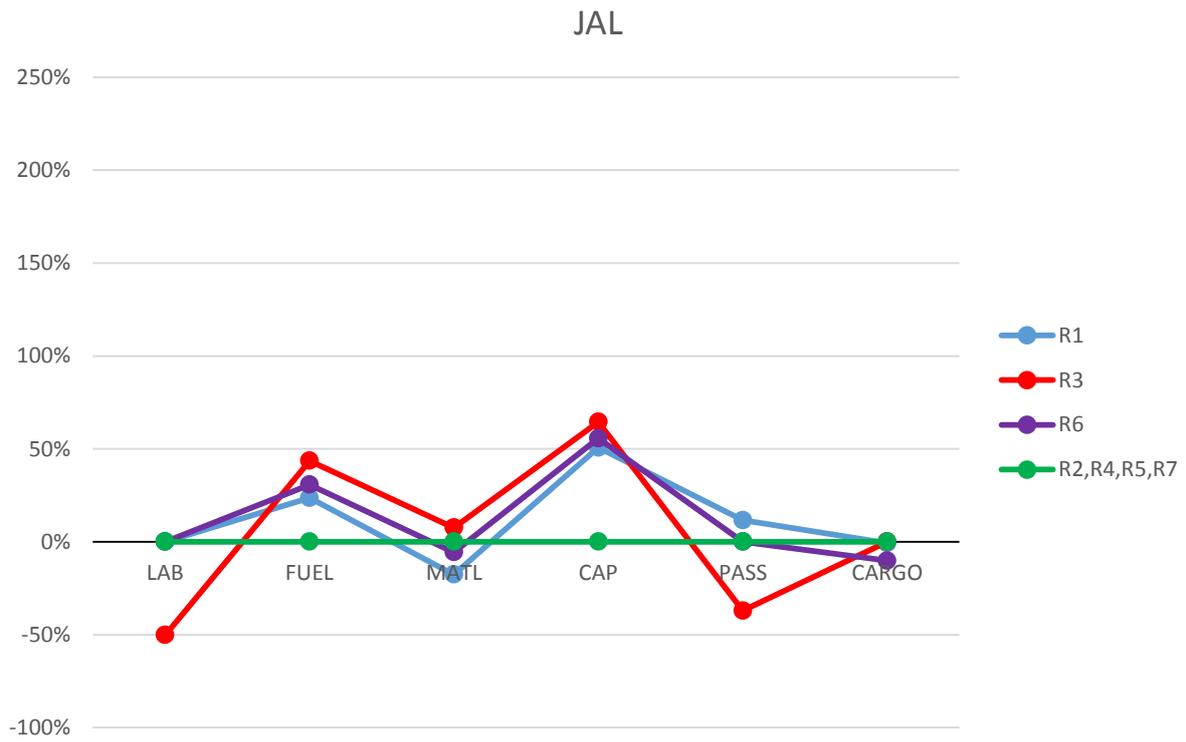

Figure 5. JAL: Deviations actual data-targets by reference set